\documentclass{amsart}
\usepackage{palatino}
\usepackage{amsthm, amsmath}
\usepackage{amssymb} 
\usepackage{amscd}
\usepackage[margin=3cm]{geometry}
\usepackage[breaklinks,bookmarksopen,bookmarksnumbered]{hyperref}
\usepackage[all]{xy}

\theoremstyle{plain}

\newtheorem{lem}{Lemma}
\newtheorem{prop}{Proposition}

\theoremstyle{remark}

\newcommand\pr{\noindent\textit{Proof} : }

\newcommand\rond{\kern 1pt{\scriptstyle\circ}\kern 1pt}

\newcommand\Ext{\operatorname{Ext}}

\newcommand\Hom{\operatorname{Hom}}

\newcommand\Ker{\operatorname{Ker}}

\newcommand\Pic{\operatorname{Pic}}
\newcommand\Div{\operatorname{Div}}

\newcommand\Z{\mathbb{Z}}

\newcommand\C{\mathbb{C}}
\renewcommand\P{\mathbb{P}}

\renewcommand\O{\mathcal{O}}

\def\qfl#1{\buildrel {#1}\over {\longrightarrow}}

\newcommand\iso{\vbox{\hbox to .8cm{\hfill{$\scriptstyle\sim$}\hfill}
\nointerlineskip\hbox to .8cm{{\hfill$\longrightarrow $\hfill}} }}

\begin{document}
\title[Vanishing thetanulls on curves with involutions]{Vanishing thetanulls on curves with involutions}
\author[Arnaud Beauville]{Arnaud Beauville}
\address{Laboratoire J.-A. Dieudonn\'e\\
UMR 7351 du CNRS\\
Universit\'e de Nice\\
Parc Valrose\\
F-06108 Nice cedex 2, France}
\email{arnaud.beauville@unice.fr}
 
\date{\today}
 
\begin{abstract}
The configuration of theta characteristics and vanishing thetanulls on a hyperelliptic curve is completely understood. We observe in this note that analogous results hold for the $ \sigma $-invariant theta characteristics on any curve $ C $ with an  involution $ \sigma $. As a consequence we get examples of  non hyperelliptic curves with a high number of vanishing thetanulls.
\end{abstract}

\maketitle 

\section{Introduction}

Let $ C $ be a smooth projective curve over $ \C $. A \emph{theta characteristic} on $ C $ is a line bundle $ \kappa $ such that $ \kappa^2\cong K_C $; it is even or odd according to the parity of $ h^0(\kappa)  $. An even theta characteristic $ \kappa $ with $ h^0(\kappa)>0  $ is called a \emph{vanishing thetanull}.

The terminology comes from the classical theory of theta functions. A theta characteristic $ \kappa $ corresponds to a symmetric theta divisor $ \Theta_\kappa $ on the Jacobian  $ JC $, defined by a theta function $ \theta_\kappa $; this function is even or odd according to the parity of $ \kappa $. Thus the numbers $ \theta_\kappa(0)  $ are 0 for $ \kappa $ odd; for $ \kappa $ even they are classical invariants attached to  the curve (``thetanullwerte" or ``thetanulls"). The number  $ \theta_\kappa(0)  $ vanishes if and only if $ \kappa $ is a vanishing thetanull in the above sense.

\par When $ C $ is hyperelliptic, the configuration of its theta characteristics and vanishing thetanulls
is completely understood (see e.g.\ \cite{M2}). We observe in this note  that analogous results  
 hold for the $ \sigma $-invariant theta characteristics on any curve $ C $ with an  involution $ \sigma $. As a consequence we obtain examples of non hyperelliptic curves with a high number of vanishing thetanulls: for instance  approximately one fourth of the even thetanulls vanish for a bielliptic curve.

 \medskip
 \section{$ \sigma $-invariant line bundles}
 
 Throughout the paper we  consider a curve $ C $ of genus $ g $, with an involution $ \sigma $. We denote by $ \pi:C\rightarrow B $ the quotient map, and by $ R\subset C $ the fixed locus of $ \sigma $. The double covering $ \pi $ determines a line bundle $ \rho $ on  $ B $ such that $ \rho^2 =\O_B(\pi_*R) $;
 we have $ \pi^*\rho=\O_C(R)  $, $ \pi_*\O_C \cong \O_B\oplus \rho^{-1} $ and $ K_C=\pi^*(K_B\otimes\rho)   $.

  For a subset $ E=\{p_1,\ldots,p_k \} $ of $ R $ we will still denote by $ E $ the divisor $ p_1+\ldots +p_k $. We consider the map $\varphi: \Z^R \rightarrow \Pic(C) $ which maps $ r\in R $ to the class of $ \O_C(r)  $. Its image lies in the subgroup $ \Pic(C)^\sigma $ of $ \sigma $-invariant line bundles.

 \begin{lem}\label{lem}
  $ \varphi $ induces a surjective homomorphism 
 $\bar\varphi: (\Z/2)^R\rightarrow \Pic(C)^\sigma/\pi^*\Pic(B)  $, whose kernel is $ \Z/2\cdot(1,\ldots,1)  $.
\end{lem}
\pr Let $ K_C $ and $ K_B $ be the fields of rational functions of $ C $ and $ B $ respectively. Let $ \langle\sigma\rangle\ (\cong\Z/2) $ be the Galois group of the covering $ \pi $. Consider the exact sequence of $ \langle\sigma\rangle $-modules
\[ 1\rightarrow K_C^*/\C^* \rightarrow \Div(C) \rightarrow \Pic(C)\rightarrow 0\ .   \]
Since $ H^1(\langle\sigma\rangle,K_C^*)=0  $ by Hilbert Theorem 90 and $ H^2(\langle\sigma\rangle,\C^*)=0  $, we have $ H^1(\langle\sigma\rangle,(K_C^*/\C^*))=0  $, hence a diagram of exact sequences:

\[  \xymatrix{1 \ar[r] & K_B^*/\C^* \ar[r]\ar[d]^\alpha & \Div(B)  \ar[r]\ar[d]^\beta & \Pic(B)  \ar[r]\ar[d]^\gamma & 0  \\
1 \ar[r] & (K_C^*/\C^*)^\sigma   \ar[r] & \Div(C)^\sigma  \ar[r]  & \Pic(C)^\sigma  \ar[r] & 0\\
} \]where the vertical arrows are induced by pull back.

If $ R= \varnothing $, this shows that $ \gamma $ is surjective, hence there is nothing to prove. Assume $ R\not= \varnothing $. Then $ \gamma $ is injective. Since $ H^1(\langle\sigma\rangle,\C^*)=\Z/2  $ and $ (K_C^*)^\sigma= K_B^*$, the cokernel of $ \alpha $ is $ \Z/2 $. The cokernel of $ \beta $ can be identified with $ (\Z/2)^R  $, so we get an exact sequence
\[ 0\rightarrow \Z/2 \longrightarrow  (\Z/2)^R \qfl{\bar\varphi} \Pic(C)^\sigma/\pi^*\Pic(B)\rightarrow 0 \ ;\]
since $ \O_C(R)\cong \pi^*\rho   $, the vector $ (1,\ldots,1)  $ belongs to $ \Ker \bar{\varphi} $, and therefore generates this kernel.\qed

\medskip
 
\begin{prop}\label{lb}
Let $ M $ be a $ \sigma $-invariant line bundle  on $ C $.

$ a) $ We have $ M \cong \pi^*L(E)  $ for some $ L\in\Pic(C) $ and $ E\subset R $. Any pair $(L',E')$ satisfying $M\cong \pi^*L'(E')$ is equal to $(L,E)$ or 
$(L\otimes  \rho^{-1}\,(\pi_*E),R-E) $. 

$ b) $ There is a natural isomorphism $ H^0(C,L)\cong  H^0(B,L)\oplus H^0(B, L\otimes  \rho^{-1}(\pi_*E))   $.
\end{prop}
\pr Part $a)$ follows directly from the lemma. Let us prove $ b) $. We view $ \O_C(E)  $ as the sheaf of rational functions on $ C $ with at most simple poles along $ E $. Then $ \sigma $ induces a homomorphism $\O_C(E)\rightarrow \sigma_*\O_C(E)  $, hence an involution of the rank 2 vector bundle $ F:=\pi_*\O_C(E) $; thus $ F $ admits a decomposition $ F=F^+\oplus F^- $ into eigen-subbundles for this involution.
 The section $ 1 $ of $ \O_C(E) $ provides a section of  $F^+ $, which generates $ F^+ $; therefore $F^-\cong \det F \cong  \rho^{-1}(\pi_*E)$. This gives a canonical decomposition $ \pi_*\O_C(E)\cong \O_B\oplus \rho^{-1}(\pi_*E) $. Taking tensor product with $ L $ and global sections gives the required isomorphism.\qed

\medskip
\section{$ \sigma $-invariant theta characteristics: the ramified case}

In this section  we assume $ R\not= \varnothing $. We denote by $ b $ the genus of $ B $ and we put $ r:= g-2b+1 $. By the Riemann-Hurwitz formula we have  $\deg \rho=r $ and $ \#R=2r $.

We now specialize Proposition \ref{lb} to the case of theta characteristics.

\begin{prop}\label{theta}
Let $ \kappa $ be a $ \sigma $-invariant theta characteristic  on $ C $.

$ a) $ We have $ \kappa \cong \pi^*L(E)  $ for some $ L\in\Pic(C) $ and $ E\subset R $ with  $ L^2\cong  K_B\otimes\rho\,(-\pi_*E) $. If another pair $(L',E')$ satisfies $\kappa \cong \pi ^*L'(E')$, we have  $(L',E')=(L,E)$ or $(L',E')=(K_B\otimes L^{-1},R-E)$.

$ b) $ We have $ h^0(\kappa)=h^0(L)+h^1(L)    $, and the parity of $ \kappa $ is equal to $ \deg(L)-(b-1 )   $. 

\end{prop}
\pr $a)$ By Proposition \ref{lb}.$ a) $ $ \kappa $ can be written $ \pi^*L(E)  $, with $ L\in\Pic(B)  $ and $ E\subset R $. The condition $ \kappa^2=K_C $ translates as  $ \pi ^*(L^2(\pi_*E))\cong \pi ^*(K_B\otimes\rho)  $. Since $R\neq\varnothing$ this implies $L^2\cong  K_B\otimes\rho\,(-\pi_*E)$. The last assertion then follows from Proposition \ref{lb}.$ a) $.

$b)$ The value of $h^0(\kappa )$ follows from Proposition \ref{lb}.$ b) $, and its parity from the Riemann-Roch theorem.\qed

\bigskip

\begin{lem}\label{lem2}
The group $ (\Pic(C)[2])^\sigma $ of $ \sigma $-invariant line bundles $ \alpha $  on $ C $ with $ \alpha^2=\O_C $ is a vector space of dimension $ 2(g-b)  $ over $ \Z/2 $.
\end{lem}
\pr  By lemma \ref{lem} we have an exact sequence
\begin{equation}\label{eq}
 0\rightarrow \Pic(B) \rightarrow \Pic(C)^\sigma \rightarrow (\Z/2)^{2r-1}\rightarrow 0\ .   
 \end{equation}
For a $ \Z $-module $ M $, let $ M[2]=\Hom(\Z/2,M)  $ be the kernel of
 the multiplication by 2 in $ M $. Note that $ \Ext^1(\Z/2,M)  $ is naturally isomorphic to $ M/2M $. Applying $ \Hom(\Z/2, -)  $ to (\ref{eq}) gives an  
exact sequence of $ (\Z/2)  $-vector spaces
\[  0\rightarrow \Pic(B)[2] \rightarrow (\Pic(C)[2])^\sigma \rightarrow (\Z/2)^{2r-1}\rightarrow \Pic(B)/2\Pic(B) \rightarrow \Pic(C)^\sigma/2\Pic(C)^\sigma \ .  \]

Let $ p\in R $. The group $ \Pic(B)/2\Pic(B) $ is generated by the class of $ \O_B(\pi p)  $; since $ \pi^*(\pi p)=2p  $, this class goes to $ 0 $ in $ \Pic(C)^\sigma/2\Pic(C)^\sigma $. Thus the dimension of 
$(\Pic(C)[2])^\sigma$ over $ \Z/2 $ is $ 2b+2r-2=\allowbreak 2(g-b)  $.\qed

\begin{prop}\label{ram}

$ a) $ The $ \sigma $-invariant theta characteristics form an affine space of dimension $ 2(g-b)  $ over $ \Z/2 $; among these, there are $ 2^{g-1}(2^{g-2b}+1) $ even theta characteristics and $ 2^{g-1}(2^{g-2b}-1) $ odd ones.

$ b) $ $ C $ admits (at least)  $ 2^{g-1}\left(2^{g-2b}+1 - 2^{-r+1}\binom{2r}{r} \right)$ vanishing thetanulls.

\end{prop}
\pr The $ \sigma $-invariant theta characteristics form  an affine space under 
$(\Pic(C)[2])^\sigma$, which has dimension $ 2(g-b)  $ by lemma \ref{lem2}.

 According to  Proposition \ref{theta}, a theta characteristic $ \kappa $ is determined by a subset $ E\subset R $  and a  line bundle $ L $ on $ B $ such that $ L^2\cong K_B\otimes\rho\,(-\pi_*E)  $. This condition  implies 
 $ \#E\equiv  r\ (\mathrm{mod.}\ 2) $.
  Moreover the parity of $ \kappa $ is that of $\deg(L)-(b-1)=   \frac{1}{2}(r-\#E)  $.

 Once $ E $ is fixed we have $ 2^{2b} $ choices for $ L $. Since  $ E $ and $ R\smallsetminus E$ give the same theta characteristic, we  consider only the subsets $ E $ with $ \#E\leq r $,  counting only half of those with $ \#E = r $. Thus the numer of even $ \sigma $-invariant theta characteristics is
 
 \[ \begin{aligned}
2^{2b} \left[\frac{1}{2} \binom{2r}{r}+\binom{2r}{r-4}+\ldots \right] &=
 2^{2b-3}\left[ (1+1)^{2r}+(-1)^r (1-1)^{2r}+(-i)^r (1+i)^{2r}+i^r(1-i)^{2r}\right]\\
& =2^{2b+2r-3}+ 2^{2b+r-2} =  2^{g-1}(2^{g-2b}+1) \ ,
 \end{aligned} \]which gives $ a) $.
 
By Proposition \ref{theta}.$ b) $ such a  theta characteristic will be a vanishing thetanull a soon as $ \deg L>b-1 $, or equivalently 
 $\#E<r$.  Thus substracting the number of  theta characteristics $ \kappa=\pi^*L(E)  $ with $\#E=r$ we obtain $ b) $.\qed
 
 \medskip	
\noindent\emph{Remarks}$.-$
1) Note that there may be more $ \sigma $-invariant vanishing thetanulls, namely those of the form $\pi ^*L(E)$ with $\deg L=b-1$ but $h^0(L)>0$. These will not occur for a general $(C,\sigma )$. 

2) Let  $g\rightarrow \infty$ with $ b $ fixed. 
By the Stirling formula $\binom{2r}{r}$ is equivalent to $2^{2r}/\sqrt{\pi r}$, so 
$  2^{-r+1}\binom{2r}{r}  $ is negligible compared to $ 2^{g-2b}=2^{r-1} $. Thus asymptotically we obtain $ 2^{2g-1-2b} $ vanishing thetanulls.

3) When $b=0$ we recover the usual numbers for hyperelliptic curves. For $b=1$ we obtain approximately $2^{2g-3}$ vanishing thetanulls, that is one fourth of the number of even theta characteristics. 

\medskip
\section{$ \sigma $-invariant theta characteristics: the \'etale case}

In this section we assume that $ \sigma $ is fixed point free ($ R=\varnothing $).

\begin{lem}
$(\Pic(C)[2])^\sigma$ is a vector space of dimension $ g+1  $ over $ \Z/2 $.
\end{lem}
\pr Apply $ \Hom(\Z/2,-)  $ to the exact sequence
\[ 0\rightarrow \Z/2 \rightarrow JB \qfl{\pi^*} JC^\sigma\rightarrow 0\ .\eqno{\qed} \]

\begin{prop} 
$ a) $ The $ \sigma $-invariant theta characteristics form an affine space of dimension $ g+1  $ over $ \Z/2 $; among these, there are $3. 2^{g-1} $ even theta characteristics and $ 2^{g-1} $ odd ones.

$ b) $ $ C $ admits a set $ \mathcal{T} $ of  $ 2^{g-2}- 2^{\frac{g-3}{2}}$ $ \sigma $-invariant vanishing thetanulls; it is contained in an affine subspace of dimension $ g-1 $ consisting of even theta characteristics.

\end{prop}
The   last property implies that for $ \kappa_1,\kappa_2,\kappa_3 $ in $ \mathcal{T} $, the theta characteristic $\kappa_1\otimes\kappa_2\otimes\kappa_3^{-1}   $ is even: 
in classical terms, $ \mathcal{T} $ is  \emph{syzygetic}. The existence of these 
vanishing thetanulls appears already in \cite{F}.
 
\medskip
\pr The first assertion follows from the previous lemma. Let $ \kappa $ be a  $ \sigma $-invariant theta characteristic; we have $ \kappa=\pi^*L $ for some line bundle $ L $ on $ C $ with $\pi^* L^2= K_C=\pi^*K_B$, which implies either $ L^2=K_B\otimes\rho $ or $ L^2=K_B $. In the first case we have
\[ h^0(\kappa)=h^0(L)+h^0(L\otimes\rho)  = h^0(L)+h^0(K_B\otimes L^{-1})\equiv 0\ \ (\mathrm{mod.}\ 2). \]
Since $\pi^* L \cong \pi^*( L\otimes\rho) $,  we get $ 2^{2b-1} $ even theta characteristics of $ C $.

\smallskip
In the second case $ L $ is a theta characteristic on $ B $. We recall briefly the theory of theta characteristics on a curve, as explained for instance in \cite{M1}. The group $ V=\Pic(B)[2]  $  is a vector space over $ \Z/2 $, equipped with a symplectic form $e$, the  \emph{Weil pairing}. A quadratic form on $ V $ associated to $ e $ is a function $ q:V\rightarrow \Z/2 $ satisfying
\[ q(\alpha+\beta)=  q(\alpha)+ q(\beta)+e(\alpha,\beta)\ . \]
The set
 $ \mathcal{Q} $ of such forms is an affine space over $ V $. 
Now the set of  theta characteristics on $ B $ is also an affine space over $ V $, which is in fact canonically 
  isomorphic to $ \mathcal{Q} $: the isomorphism associates 
 to a theta characteristic $ L $  the form $ q_L\in \mathcal{Q} $ defined by $ q_L(\alpha)=h^0(L\otimes\alpha)+h^0(L) \ \ (\mathrm{mod.}\ 2)$. Moreover the parity of $ L $ is equal to the Arf invariant  $\mathrm{Arf}( q_L) $.

\smallskip
Coming back to our situation, let $ L $ be a theta characteristic  on $ B $, and $ \kappa=\pi^*L $; we have
\[ h^0(\kappa) = h^0(L)+h^0(L\otimes\rho) \equiv q_L(\rho)\  \ (\mathrm{mod.}\ 2).   \]
The function $ q\mapsto q(\rho)  $ is an affine function on $\mathcal{Q}  $, hence it takes equally often the values $ 0 $ and $ 1 $. Taking into account the isomorphism $ \pi^*L\cong \pi^*(L\otimes\rho) $, we get $ 2^{2b-2} $ even 
theta characteristics on $ C $ and $ 2^{2b-2} $ odd ones; summing up we obtain $ a) $.

\smallskip
Suppose $ \kappa =\pi^*L $ is even, that is, $ h^0(L)\equiv h^0(L\otimes\rho )\   (\mathrm{mod.}\ 2)$ ; if we want $ h^0(\kappa)>0  $, a good way (actually the only one if $ B $ is generic) is to choose $ L $ odd, that is, $\mathrm{Arf}( q_L) =1$. Equivalently, we look for forms $ q\in \mathcal{Q} $ with $ q(\rho)=0  $ and $ \mathrm{Arf}( q)=1 $.

Let $ \rho' $ be an element of $V $ with $ e(\rho,\rho')=1  $. $ \rho $ and $ \rho' $ span a plane $ P \subset V$, such that $ V=P\oplus P^\perp $. A form $ q\in \mathcal{Q} $ is determined by its restriction to $ P $ and $ P^\perp $, and we have $ \mathrm{Arf}( q)=\mathrm{Arf}( q_{|P})+\mathrm{Arf}( q_{|P^\perp}) $. The condition  $ q(\rho) =0 $ implies $ \mathrm{Arf}( q_{|P})=q(\rho)q'(\rho)=  0 $; so $ q $ is determined by $ q(\rho')\in \Z/2 $ and a form $ q' $ on $ P^\perp $ with Arf invariant 1. Since $ \dim P^\perp=2(b-1)  $, there are $ 2^{b-2}(2^{b-1}-1) $ such forms, hence $ 2^{b-1}(2^{b-1}-1) $ forms $ q\in \mathcal{Q} $ with $ q(\rho)=0  $ and $ \mathrm{Arf}( q)=1 $. Taking again into account the isomorphism $ \pi^*L\cong \pi^*(L\otimes\rho) $, we obtain $ 2^{b-2}(2^{b-1}-1) = 2^{g-2}- 2^{\frac{g-3}{2}}$ vanishing thetanulls on $ C $. 

They are contained in the affine space of  theta characteristics $ \kappa=\pi^*L $ with $ q_L{}(\rho)=0  $, which has dimension $ 2b-2=g-1 $ and consists of even
theta characteristics.\qed

 \medskip
 \section{Low genus}
Let $ C $ be a non hyperelliptic curve of genus $ g $. How many  vanishing thetanulls can $ C $ have? 
The answer is well-known up to genus 5. There is no vanishing thetanull in genus 3, and at most one in genus 4 (which occurs if and only if the unique quadric containing the canonical curve is singular).  

Suppose $ g=5 $. If $ C $ is trigonal it admits at most one vanishing thetanull. Otherwise the canonical curve $ C\subset \P^4 $ is the base locus of a net $ \Pi $ of quadrics. The discriminant curve  (locus of the quadrics in $ \Pi $ of rank $ \leq 4 $) is a plane quintic with only ordinary nodes; these nodes  correspond to the rank 3 quadrics of $ \Pi $, that is to the vanishing thetanulls of $ C $. Therefore $ C $ can have any number $ \leq 10 $ of vanishing thetanulls; they are syzygetic \cite{A}. The maximum 10 is attained by the so-called Humbert curves, for which all the quadrics in $ \Pi $ can be simultaneously diagonalized. 
They have an action of the group $ (\Z/2)^4  $, generated by 5 involutions with elliptic quotient.

\smallskip
Starting with $ g=6 $ very little seems to be known.  By Proposition \ref{ram}.$ b) $, if $ C $ is bielliptic (that is, $ C $ admits an involution with elliptic quotient), it has 40 vanishing thetanulls. This can be slightly improved  as follows. We take an elliptic curve $ B $, a line bundle $ \alpha $ of degree 2  on $ B $, a point $ p \in B$, and disjoint divisors $ A $ in $ |\alpha(p) | $, $ A_1,A_2,A_3 $ in $ |\alpha| $  which do not contain $ p $. We put $ \rho=\alpha^2(p)  $ and $ \bar R=A_1+A_2+A_3+A+p $, and construct the double covering $ \pi:C\rightarrow B $ associated to $ (\rho,\bar R)  $. 
  The curve $ C $ has three extra vanishing thetanulls, namely $ \O_C(\tilde A_i+\tilde{A}_j+\tilde p)  $ for $  i<j $, where  $ \tilde{A}_i$ and $ \tilde{p} $ are the lifts of $ A_i$ and $ p $ to $ C $. Thus we get a genus 6 curve with 43 vanishing thetanulls;  it is likely that one can do better.

\end{document}